\documentclass[a4paper,11pt]{amsart}
\usepackage{graphicx}
\usepackage{mathptmx}
\usepackage{amsmath}
\usepackage{amssymb}
\usepackage{enumitem}
\usepackage{xcolor}
\usepackage{xparse}
\NewDocumentCommand{\eulerian}{omm}
 {%
  \genfrac<>{0pt}{}{#2}{#3}%
  \IfValueT{#1}{_{\!#1}}%
 }

 \newmuskip\pFqmuskip

\newcommand*\pFq[6][8]{%
  \begingroup % only local assignments
  \pFqmuskip=#1mu\relax
  \mathchardef\normalcomma=\mathcode`,
  \mathcode`\,=\string"8000
  \begingroup\lccode`\~=`\,
  \lowercase{\endgroup\let~}\pFqcomma
  {}_{#2}F_{#3}{\left(\genfrac..{0pt}{}{#4}{#5}\bigg|#6\right)}%
  \endgroup
}
\newcommand{\pFqcomma}{{\normalcomma}\mskip\pFqmuskip}

\newtheorem{theorem}{Theorem}

\newtheorem{corollary}[theorem]{Corollary}
\newtheorem{proposition}[theorem]{Proposition}

\begin{document}

\title[Some results on degenerate Fubini and degenerate Bell polynomials]{Some results on degenerate Fubini and degenerate Bell polynomials}

\author{Taekyun  Kim}
\address{Department of Mathematics, Kwangwoon University, Seoul 139-701, Republic of Korea}
\email{tkkim@kw.ac.kr}

\author{DAE SAN KIM}
\address{Department of Mathematics, Sogang University, Seoul 121-742, Republic of Korea}
\email{dskim@sogang.ac.kr}


\maketitle

\begin{abstract}
The aim of this paper is to further study some properties and identities on the degenerate Fubini and the degenerate Bell polynomials which are degenerate versions of the Fubini and the Bell polynomials, respectively. Especially, we find several expressions for the generating function of the sum of the values of the generalized falling factorials at positive consecutive integers.

\end{abstract}

\section{Introduction}

Carlitz [3] initiated a study of degenerate versions of Bernoulli and Euler polynomials,
namely the degenerate Bernoulli and the degenerate Euler polynomials. In recent years, some mathematicians have explored degenerate versions of many special polynomials and numbers, which include the degenerate Bernoulli numbers of the second kind, the degenerate Stirling numbers of the first and second kinds, the degenerate Bell numbers and polynomials, and so on (see [6,8,10-12,14-17,19] and the references therein).
Many interesting results on such degenerate versions together with their applications were found by employing several different tools such as combinatorial methods, generating functions, $p$-adic analysis, umbral calculus techniques, differential equations, probability theory, operator theory, special functions and analytic number theory. It is remarkable that the exploration for degenerate versions is not only limited to polynomials and numbers but also extended to transcendental functions like gamma functions (see [13]).\par
Recently, the degenerate Fubini polynomials and the degenerate Bell polynomials are introduced (see [12,16]).
The aim of this paper is to further study some properties and identities on the degenerate Fubini and the degenerate Bell polynomials which are degenerate versions of the Fubini and the Bell polynomials, respectively. Especially, we find several expressions for the generating functions of the sum $S_{n,\lambda}(p)=(1)_{p,\lambda}+(2)_{p,\lambda}+\cdots+(n)_{p,\lambda}$ of the values of the generalized falling factorials at positive consecutive integers.
The outline of this paper is as follows. In Section 1, we recall degenerate exponential functions, degenerate logarithms, degenerate Stirling numbers of the first kind and degenerate Stirling numbers of the second kind which are respectively degenerate versions of exponentials, logarithms, Stirling numbers of the first kind and Stirling numbers of the second kind. Also, we recall the degenerate Bell polynomials, the degenerate polylogarithm and the degenerate Fubini polynomials.
In Section 2, we derive recurrence relations for the degenerate Bell polynomilas in Theorems 1 and 2. In Theorem 4 and 7, the exponential generating function of the sum $S_{n,\lambda}(p)$ is determined in terms of the degenerate Bell polynomials. Another exponential generating function for the sum $S_{n,\lambda}(p)$ is obtained in Theorem 6 in terms of the degenerate Stirling numbers of the second kind. In Theorem 8, the generating function for $S_{n,\lambda}(p)$ is determined in terms of the degenerate Fubini polynomial. In Theorem 12, obtained is a differential equation satisfied by the degenerate Fubini polynomials. Finally, we conclude our paper in Section 3. \par

The generalized falling factorial sequence is defined by
\begin{equation}
(x)_{0,\lambda}=1,\quad (x)_{n,\lambda}=x(x-\lambda)(x-2\lambda)\cdots\big(x-(n-1)\lambda\big),\quad (n\ge 1),\label{1}
\end{equation}
where $\lambda$ is a real number (see [8,12]). \par
For any $\lambda\in\mathbb{R}$, the degenerate exponentials are defined by
\begin{equation}
e_{\lambda}^{x}(t)=\sum_{n=0}^{\infty}(x)_{n,\lambda}\frac{t^{n}}{n!},\quad (\mathrm{see}\ [11,14]).\label{2}
\end{equation}
In particular, for $x=1$, $\displaystyle e_{\lambda}(t)=e_{\lambda}^{1}(t)=\sum_{n=0}^{\infty}(1)_{n,\lambda}\frac{t^{n}}{n!}\displaystyle$. Note that $\displaystyle\lim_{\lambda\rightarrow 0}e_{\lambda}^{x}(t)=e^{xt}\displaystyle$. \par
The degenerate Bell polynomials are considered by Kim-Kim and given by
\begin{equation}
e^{x(e_{\lambda}(t)-1)}=\sum_{n=0}^{\infty}\phi_{n,\lambda}(x)\frac{t^{n}}{n!},\quad (\mathrm{see}\ [14,15]).\label{3}
\end{equation}
From \eqref{3}, we note that
\begin{displaymath}
	\lim_{n\rightarrow\infty}\phi_{n,\lambda}(x)=\phi_{n}(x),\quad (n\ge 0),
\end{displaymath}
where $\phi_{n}(x)$ are the ordinary Bell polynomials given by
\begin{equation}
e^{x(e^{t}-1)}=\sum_{n=0}^{\infty}\phi_{n}(x)\frac{t^{n}}{n!},\quad (\mathrm{see}\ [6-12,14-20]).\label{4}
\end{equation}
In [8], the degenerate Stirling numbers of the first kind are defined by
\begin{equation}
(x)_{n}=\sum_{k=0}^{n}S_{1,\lambda}(n,k)(x)_{k,\lambda},\quad (n\ge 0),\label{5}	
\end{equation}
where $(x)_{0}=1,\ (x)_{n}=x(x-1)\cdots(x-n+1),\ (n\ge 1)$. \par
Note that $\displaystyle \lim_{\lambda\rightarrow 9}S_{1,\lambda}(n,k)=S_{1}(n,k)\displaystyle$ are the ordinary Stirling number of the first kind given by
\begin{displaymath}
(x)_{n}=\sum_{k=0}^{n}S_{1}(n,k)x^{k},\quad (n\ge 0),\quad (\mathrm{see}\ [1-5,7]).
\end{displaymath}
As the inversion formula of \eqref{5}, the degenerate Stirling numbers of the second kind are introduced by Kim-Kim and given by
\begin{equation}
(x)_{n,\lambda}=\sum_{k=0}^{n}S_{2,\lambda}(n,k)(x)_{k},\quad (n\ge 0),\quad (\mathrm{see}\ [8]).\label{6}
\end{equation}
Let $\log_{\lambda}t$ be the compositional inverse of $e_{\lambda}(t)$. Then we have
\begin{equation}
	\log_{\lambda}(1+t)=\sum_{n=1}^{\infty}\frac{\lambda^{n-1}(1)_{n,1/\lambda}}{n!}t^{n},\quad (\mathrm{see}\ [8]).\label{7}
\end{equation}
Note that $\displaystyle\lim_{\lambda\rightarrow 0}\log_{\lambda}(1+t)=\log(1+t)\displaystyle$. \par
For $k\in\mathbb{Z}$, the degenerate polylogarithm function is defined by Kim-Kim and given by
\begin{equation}
\mathrm{Li}_{k,\lambda}(t)=\sum_{n=1}^{\infty}\frac{(-\lambda)^{n-1}(1)_{n,1/\lambda}}{(n-1)!n^{k}}t^{n},\quad (\mathrm{see}\ [8]).\label{8}	
\end{equation}
Note that $\displaystyle\lim_{\lambda\rightarrow 0}\mathrm{Li}_{k,\lambda}(t)=\mathrm{Li}_{k}(t)=\sum_{n=1}^{\infty}\frac{x^{n}}{n^{k}}\displaystyle$. is the polylogarithm function and $\mathrm{Li}_{1,\lambda}(t)=-\log_{\lambda}(1-t)$. \par
From \eqref{5} and \eqref{10}, we can derive the generating functions of degenerate Stirling numbers which are given by
\begin{equation}
\frac{1}{k!}\big(\log_{\lambda}(1+t)\big)^{k}=\sum_{n=k}^{\infty}S_{1,\lambda}(n,k)\frac{t^{n}}{n!},\label{9}
\end{equation}
and
\begin{equation}
\frac{1}{k!}\big(e_{\lambda}(t)-1\big)^{k}=\sum_{n=k}^{\infty}S_{2,\lambda}(n,k)\frac{t^{n}}{n!},\quad (k\ge 0),\quad (\mathrm{see}\ [8]).\label{10}
\end{equation}
In [15], the degenerate Fubini polynomials are defined by
\begin{equation}
\frac{1}{1-x(e_{\lambda}(t)-1)}=\sum_{n=0}^{\infty}F_{n,\lambda}(x)\frac{t^{n}}{n!}.\label{11}
\end{equation}
Thus, by \eqref{11}, we get
\begin{equation}
F_{n,\lambda}(x)=\sum_{k=0}^{n}k!S_{2,\lambda}(n,k)x^{k},\quad (n\ge 0),\quad (\mathrm{see}\ [15]).\label{12}
\end{equation}
By \eqref{3}, we easily get
\begin{align}
\phi_{n,\lambda}(x)&=\sum_{k=0}^{n}S_{2,\lambda}(n,k)x^{k},\quad (n\ge 0), \label{13} \\
&=e^{-x}\sum_{k=0}^{\infty}\frac{(k)_{n,\lambda}}{k!}x^{k},\quad (\mathrm{see}\ [9,10,14,15]).\nonumber
\end{align}

\section{Some new formulae on degenerate Fubini and Bell polynomials}
First, by \eqref{3}, we get
\begin{align}
\frac{d}{dx}e^{x(e_{\lambda}(t)-1)}&=(e_{\lambda}(t)-1)e^{x(e_{\lambda}(t)-1)} \label{14} \\
&=\sum_{n=0}^{\infty}\bigg(\sum_{k=0}^{n}\binom{n}{k}(1)_{n-k,\lambda}\phi_{k,\lambda}(x)-\phi_{n,\lambda}(x)\bigg)\frac{t^{n}}{n!}.\nonumber	
\end{align}
Thus, by \eqref{3} and \eqref{14}, we get
\begin{align}
\phi^{\prime}_{n,\lambda}(x)&=\frac{d}{dx}\phi_{n,\lambda}(x)=\sum_{k=0}^{n}\binom{n}{k}(1)_{n-k,\lambda}\phi_{k,\lambda}(x)-\phi_{n,\lambda}(x) \label{15}\\
&=\sum_{k=0}^{n-1}\binom{n}{k}(1)_{n-k,\lambda}\phi_{k,\lambda}(x),\quad (n\ge 1). \nonumber
\end{align}
From \eqref{3}, we note that
\begin{align}
\sum_{n=0}^{\infty}\phi_{n+1,\lambda}(x)\frac{t^{n}}{n!}&=\frac{d}{dt}e^{x(e_{\lambda}(t)-1)}=xe_{\lambda}^{1-\lambda}(t)e^{x(e_{\lambda}(t)-1)}\label{16} \\
&=x \sum_{l=0}^{\infty}(1-\lambda)_{l,\lambda}\frac{t^{l}}{l!}\sum_{k=0}^{\infty}\phi_{k,\lambda}(x)\frac{t^{k}}{k!} \nonumber \\
&=\sum_{n=0}^{\infty}\bigg(x\sum_{k=0}^{n}\binom{n}{k}\phi_{k,\lambda}(x)(1-\lambda)_{n-k,\lambda}\bigg)\frac{t^{n}}{n!}. \nonumber
\end{align}
By comparing the coefficients on both sides of \eqref{15}, we have
\begin{equation}
\phi_{n+1,\lambda}(x)=x\sum_{k=0}^{n}\binom{n}{k}\phi_{k,\lambda}(x)(1-\lambda)_{n-k,\lambda},\quad (n\ge 0).\label{17}
\end{equation}
From the definition in \eqref{1}, we easily get
\begin{align}
(1-\lambda)_{n-k,\lambda}=(1)_{n-k,\lambda}\big(1-(n-k)\lambda\big)\label{18}.
\end{align}
From \eqref{17} and \eqref{18}, we have
\begin{align}
\phi_{n+1,\lambda}(x)&= x\sum_{k=0}^{n}\binom{n}{k}\phi_{k,\lambda}(x)(1)_{n-k,\lambda}\big(1-(n-k)\lambda\big)\label{19} \\
&= x\sum_{k=0}^{n}\binom{n}{k}\phi_{k,\lambda}(x)(1)_{n-k,\lambda}-x\lambda\sum_{k=0}^{n}\binom{n}{k}(n-k)(1)_{n-k,\lambda}\phi_{k,\lambda}(x)\nonumber \\
&= x\sum_{k=0}^{n}\binom{n}{k}\phi_{k,\lambda}(x)(1)_{n-k,\lambda}-nx\lambda\sum_{k=0}^{n}\binom{n-1}{k}\phi_{k,\lambda}(x)(1)_{n-k,\lambda}.\nonumber
\end{align}
Therefore, by \eqref{19}, we obtain the following theorem.
\begin{theorem}
For $n\in\mathbb{Z}$ with $n\ge 0$, we have 	\begin{displaymath}
	\sum_{k=0}^{n}\binom{n}{k}\phi_{n,k}(x)(1)_{n-k,\lambda}=\frac{1}{x}\phi_{n+1,\lambda}(x)+n\lambda\sum_{k=0}^{n-1}\binom{n-1}{k}\phi_{k,\lambda}(x)(1)_{n-k,\lambda}.
\end{displaymath}
\end{theorem}
By \eqref{15} and Theorem 1, we obtain the following theorem.
\begin{theorem}
For $n\in\mathbb{Z}$ with $n\ge 0$, we have
	\begin{displaymath}
		\frac{1}{x}\phi_{n+1,\lambda}(x)=
\phi_{n,\lambda}^{\prime}(x)+\phi_{n,\lambda}(x)-n\lambda\sum_{k=0}^{n-1}\binom{n-1}{k}\phi_{k,\lambda}(x)(1)_{n-k,\lambda}.
	\end{displaymath}
\end{theorem}
Let
\begin{equation}
\begin{aligned}
	Y_{\lambda}(x,p)&=\sum_{n=1}^{\infty}\sum_{k=1}^{n}(k)_{p,\lambda}\frac{x^{n}}{n!} \\
	&=\sum_{n=1}^{\infty}\Big((1)_{p,\lambda}+(2)_{p,\lambda}+\cdots+(n)_{p,\lambda}\Big)\frac{x^{n}}{n!}
\end{aligned}	\label{21}
\end{equation}
where $p$ is a positive integer. \par
From \eqref{21}, we note that
\begin{align}
Y_{\lambda}(x,p)=\sum_{k=1}^{\infty}(k)_{p,\lambda}\sum_{n=k}^{\infty}\frac{x^{n}}{n!}=\sum_{k=1}^{\infty}(k)_{p,\lambda}\bigg(e^{x}-\sum_{l=0}^{k-1}\frac{x^{l}}{l!}\bigg). \label{22}
\end{align}
In \eqref{13}, we note that
\begin{equation}
e^{x}\phi_{p,\lambda}(x)=\sum_{n=1}^{\infty}\frac{(n)_{p,\lambda}}{n!}x^{n},\quad (p\in\mathbb{Z}).
\end{equation}
Let
\begin{equation}
\begin{aligned}
	y_{p,\lambda}(x)&=Y_{\lambda}(x,p)-e^{x}\phi_{p,\lambda}(x)\\
	&=\sum_{n=2}^{\infty}\Big((1)_{p,\lambda}+(2)_{p,\lambda}+\cdots+(n-1)_{p,\lambda}\Big)\frac{x^{n}}{n!},\quad (p\in\mathbb{N}).
\end{aligned}	\label{24}
\end{equation}
By \eqref{24}, we get
\begin{align}
y_{p,\lambda}^{\prime}(x)=\frac{d}{dx}y_{p,\lambda}(x)&=\sum_{n=2}^{\infty}\Big((1)_{p,\lambda}+(2)_{p,\lambda}+\cdots+(n-1)_{p,\lambda}\Big)\frac{x^{n-1}}{(n-1)!}\label{25}\\
&=\sum_{n=1}^{\infty}\Big((1)_{p,\lambda}+(2)_{p,\lambda}+\cdots+(n)_{p,\lambda}\Big)\frac{x^{n}}{n!}\nonumber \\
&=\sum_{n=1}^{\infty}\Big((1)_{p,\lambda}+\cdots+(n-1)_{p,\lambda}\Big)\frac{x^{n}}{n!}+\sum_{n=1}^{\infty}(n)_{p,\lambda}\frac{x^{n}}{n!}\nonumber \\
&=y_{p,\lambda}(x)+\sum_{n=1}^{\infty}(n)_{p,\lambda}\frac{x^{n}}{n!}=y_{p,\lambda}(x)+e^{x}\phi_{p,\lambda}(x).\nonumber
\end{align}
From \eqref{25}, we have the following differential equation.
\begin{proposition}
For $p\in\mathbb{N}$, we have
\begin{displaymath}
y_{p,\lambda}^{\prime}(x)-y_{p,\lambda}(x)=e^{x}\phi_{p,\lambda}(x).
\end{displaymath}
\end{proposition}
By Proposition 3, we easily get
\begin{equation}
\frac{d}{dx}\Big(e^{-x}y_{p,\lambda}(x)\Big)=\phi_{p,\lambda}(x).\label{27}	
\end{equation}
Thus, by \eqref{27}, we get
\begin{equation}
y_{p,\lambda}(x)=e^{x}\int_{0}^{x}\phi_{p,\lambda}(t)dt.\label{29}	
\end{equation}
Therefore, by \eqref{21}, \eqref{24} and \eqref{29}, we obtain the following theorem.
\begin{theorem}
For $p\in\mathbb{N}$, we have
\begin{displaymath}
\sum_{n=1}^{\infty}\Big((1)_{p,\lambda}+ (2)_{p,\lambda}+\cdots+ (n)_{p,\lambda}\Big)\frac{t^{n}}{n!}=e^{x}\phi_{p,\lambda}(x)+e^{x}\int_{0}^{x}\phi_{p,\lambda}(t)dt.
\end{displaymath}	
\end{theorem}
In [10], considered are the degenerate differential operators which are given by
\begin{equation}
\Big(x\frac{d}{dx}\Big)_{n,\lambda}=\Big(x\frac{d}{dx}\Big)\Big(x\frac{d}{dx}-\lambda\Big)\cdots\Big(x\frac{d}{dx}-(n-1)\lambda\Big),\quad (n\ge 1).\label{30}	
\end{equation}
For $p\in\mathbb{N}$, by \eqref{30}, we get
\begin{align}
\Big(x\frac{d}{dx}\Big)_{p,\lambda}e^{x}=\sum_{n=1}^{\infty}\frac{(n)_{p,\lambda}}{n!}x^{n}=e^{x}\phi_{p,\lambda}(x). \label{31}
\end{align}
Therefore, by \eqref{31}, we obtain the following theorem.
\begin{theorem}
For $p\in\mathbb{N}$, we have
\begin{displaymath}
\Big(x\frac{d}{dx}\Big)_{p,\lambda}e^{x}=\sum_{n=1}^{\infty}\frac{(n)_{p,\lambda}}{n!}x^{n}=e^{x}\phi_{p,\lambda}(x).
\end{displaymath}	
\end{theorem}
From \eqref{24} and \eqref{29}, we note that
\begin{align}
Y_{\lambda}(x,p)&=e^{x}\phi_{p,\lambda}(x)+e^{x}\int_{0}^{x}\phi_{p,\lambda}(t)dt 	\label{32} \\
&=e^{x}\phi_{p,\lambda}(x)+e^{x}\sum_{k=0}^{p}S_{2,\lambda}(p,k)\int_{0}^{x}t^{k}dt \nonumber \\
&=e^{x}\phi_{p,\lambda}(x)+e^{x}\sum_{k=0}^{p}S_{2,\lambda}(p,k)\frac{x^{k+1}}{k+1}\nonumber \\
&=e^{x}\sum_{k=0}^{p}\bigg(x{^k}+\frac{x^{k+1}}{k+1}\bigg)S_{2,\lambda}(p,k).\nonumber
\end{align}
Therefore, by \eqref{21} and \eqref{32}, we obtain the following theorem.
\begin{theorem}
For $p\in\mathbb{N}$, we have
\begin{displaymath}
	\sum_{n=1}^{\infty}\Big((1)_{p,\lambda}+(2)_{p,\lambda}+\cdots+(n)_{p,\lambda}\Big)\frac{x^{n}}{n!}=e^{x}\sum_{k=0}^{p}\bigg(x^{k}+\frac{x^{k+1}}{k+1}\bigg)S_{2,\lambda}(p,k).
\end{displaymath}	
\end{theorem}
By \eqref{6}, we easily get
\begin{equation}
S_{2,\lambda}(n+1,k)=S_{2,\lambda}(n,k-1)+(k-n\lambda)S_{2,\lambda}(n,k),\label{33}
\end{equation}
where $n,k\ge 0$ with $n\ge k$. \par
From \eqref{31}, \eqref{30} and noting $S_{2,\lambda}(p,k)=0$ for $p \ge 1$, we note that
\begin{align}
Y_{\lambda}(x,p)&=e^{x}\sum_{k=0}^{p}\bigg(x^{k}+\frac{x^{k+1}}{k+1}\bigg)S_{2,\lambda}(p,k)\label{34} \\
&=e^{x}\sum_{k=1}^{p+1}x^{k}S_{2,\lambda}(p,k)+e^{x}\sum_{k=0}^{p}\frac{x^{k+1}}{k+1}S_{2,\lambda}(p,k)\nonumber \\
&=e^{x}\sum_{k=1}^{p+1}\bigg(S_{2,\lambda}(p,k)(k-p\lambda)+S_{2,\lambda}(p,k-1)\bigg)\frac{x^{k}}{k}+p\lambda e^{x}\sum_{k=1}^{p+1}\frac{S_{2,\lambda}(p,k)}{k}x^{k}\nonumber \\
&=e^{x}\sum_{k=1}^{p+1}\frac{1}{k}S_{2,\lambda}(p+1,k)x^{k}+p\lambda e^{x}\sum_{k=1}^{p+1}\frac{1}{k}S_{2,\lambda}(p,k)x^{k}\nonumber \\
&=e^{x}\sum_{k=1}^{p+1}\bigg(\frac{S_{2,\lambda}(p+1,k)}{k}+\frac{p \lambda }{k}S_{2,\lambda}(p,k)\bigg)x^{k}\nonumber \\
&=e^{x}\int_{0}^{x}\Big(\phi_{p+1,\lambda}(t)+p\lambda \phi_{p,\lambda}(t)\Big)\frac{dt}{t},\nonumber
\end{align}
where $p\in\mathbb{N}$ and $\displaystyle\frac{1}{t}\phi_{p+1,\lambda}(t)=\sum_{k=1}^{p+1}S_{2,\lambda}(p+1,k)t^{k-1}\displaystyle$. \par
Therefore, by \eqref{34}, we obtain the following theorem.
\begin{theorem}
For $p\in\mathbb{Z}$, we have 	
\begin{displaymath}
	\sum_{n=1}^{\infty}\Big((1)_{p,\lambda}+(2)_{p,\lambda}+\cdots+(n)_{p,\lambda}\Big)\frac{x^{n}}{n!}= e^{x}\int_{0}^{x}\Big(\phi_{p+1,\lambda}(t)+p\lambda \phi_{p,\lambda}(t)\Big)\frac{dt}{t}.
\end{displaymath}
\end{theorem}
From \eqref{12}, we note that
\begin{align}
\sum_{l=0}^{\infty}(l)_{p,\lambda}x^{l}&=\sum_{l=0}^{\infty}x^{l}\sum_{k=0}^{p}S_{2,\lambda}(p,k)(l)_{k}= \sum_{l=0}^{\infty}x^{l}\sum_{k=0}^{p}S_{2,\lambda}(p,k)\binom{l}{k}k!\label{35}\\
&=\sum_{l=k}^{\infty}x^{l}\sum_{k=0}^{p}S_{2,\lambda}(p,k)\binom{l}{k}k!=\sum_{k=0}^{p}S_{2,\lambda}(p,k)k!\sum_{l=0}^{\infty}\binom{l+k}{k}x^{l+k}\nonumber \\
&=\sum_{k=0}^{p}S_{2,\lambda}(p,k)k!x^{k}\sum_{l=0}^{\infty}\binom{l+k}{k}x^{l}=\sum_{k=0}^{p}S_{2,\lambda}(p,k)k!x^{k}\Big(\frac{1}{1-x}\Big)^{k+1}\nonumber \\
&=\frac{1}{1-x}\sum_{k=0}^{p}S_{2,\lambda}(p,k)k!\Big(\frac{x}{1-x}\Big)^{k}=\frac{1}{1-x}F_{p,\lambda}\Big(\frac{x}{1-x}\Big). \nonumber
\end{align}
For $p\in\mathbb{Z}$ with $p\ge 0$, and by using \eqref{35}, we have
\begin{equation}
\Big(x\frac{d}{dx}\Big)_{p,\lambda}\frac{1}{1-x}=\sum_{n=0}^{\infty}\Big(x\frac{d}{dx}\Big)_{p,\lambda}x^{n}=\sum_{n=0}^{\infty}(n)_{p,\lambda}x^{n}=\frac{1}{1-x}F_{p,\lambda}\Big(\frac{x}{1-x}\Big).\label{36}
\end{equation}
Thus, for $p\in\mathbb{N}$, we have
\begin{align}
&\sum_{n=1}^{\infty}\Big((1)_{p,\lambda}+(2)_{p,\lambda}+\cdots+(n)_{p,\lambda}\Big)x^{n}=\sum_{n=1}^{\infty}x^{n}\sum_{k=1}^{n}(k)_{p,\lambda}\label{37} \\
&=\sum_{k=1}^{\infty}(k)_{p,\lambda}\sum_{n=k}^{\infty}x^{n}=\frac{1}{1-x}\sum_{k=1}^{\infty}(k)_{p,\lambda}x^{k}=\frac{1}{(1-x)^{2}}F_{p,\lambda}\Big(\frac{x}{1-x}\Big).\nonumber
\end{align}
Therefore, by \eqref{36} and \eqref{37}, we obtain the following theorem.
\begin{theorem}
For $p\ge 0$, we have
\begin{displaymath}
\Big(x\frac{d}{dx}\Big)_{p,\lambda}\frac{1}{1-x}=\sum_{n=0}^{\infty}(n)_{p,\lambda}x^{n}=\frac{1}{1-x}F_{p,\lambda}\Big(\frac{x}{1-x}\Big).
\end{displaymath}	
In addition, for $p\in\mathbb{N}$, we have
\begin{displaymath}
\sum_{n=1}^{\infty}\Big((1)_{p,\lambda}+(2)_{p,\lambda}+\cdots+(n)_{p,\lambda}\Big)x^{n}=\frac{1}{(1-x)^{2}}F_{p,\lambda}\Big(\frac{x}{1-x}\Big).
\end{displaymath}
\end{theorem}
Now, we observe that, for $p\in\mathbb{N}$,
\begin{align}
\sum_{n=1}^{\infty}(n)_{p,\lambda}&\bigg(\frac{1}{1-x}-\sum_{l=0}^{n}x^{l}\bigg)
=\sum_{n=1}^{\infty}(n)_{p,\lambda}\frac{x^{n+1}}{1-x}\label{38} \\
&=\frac{x}{1-x}\sum_{n=1}^{\infty}(n)_{p,\lambda}x^{n}
=\frac{x}{(1-x)^2}F_{p,\lambda}\Big(\frac{x}{1-x}\Big).\nonumber
\end{align}
Therefore, by \eqref{38}, we obtain the following corollary.
\begin{corollary}
For $p\in\mathbb{N}$, we have
\begin{displaymath}
\sum_{n=1}^{\infty}(n)_{p,\lambda}\bigg(\frac{1}{1-x}-\sum_{l=0}^{n}x^{l}\bigg)= \frac{x}{(1-x)^2}F_{p,\lambda}\Big(\frac{x}{1-x}\Big).
\end{displaymath}
\end{corollary}
From \eqref{8}, we note that
\begin{align}
\frac{1}{1-x}\mathrm{Li}_{p,\lambda}(x)&=\sum_{l=0}^{\infty}x^{l}\sum_{k=1}^{\infty}\frac{(-\lambda)^{k-1}(1)_{k,1/\lambda}}{(k-1)!k^{p}}x^{k}\label{39} \\
&=\sum_{n=1}^{\infty}\sum_{k=1}^{n}\frac{(-\lambda)^{k-1}(1)_{k,1/\lambda}}{k^{p}(k-1)!}x^{n}.\nonumber	
\end{align}
For any nonnegative integer $r$, let $D_{r}$ denote the linear operator $D_{r}=\frac{1}{r!}\big(\frac{d}{dx}\big)^{r}x^{r}$, so that, for any polynomial $p(x) \in \mathbb{C}[x]$, we have
\begin{align}
D_{r}p(x)=\frac{1}{r!}\Big(\frac{d}{dx}\Big)^{r}\bigg[x^{r}p(x)\bigg].\label{40}
\end{align}
From \eqref{39}, we note that
\begin{align}
D_{r}F_{n,\lambda}(x)&=\frac{1}{r!}\bigg(\frac{d}{dx}\bigg)^{r}\bigg[\sum_{k=0}^{n}S_{2,\lambda}(n,k)k!x^{k+r}\bigg] \label{41} \\
&=\frac{1}{r!}\bigg[\sum_{k=0}^{n}S_{2,\lambda}(n,k)k!(k+r)_{r}x^{k}\bigg]\nonumber \\
&=\sum_{k=0}^{n}S_{2,\lambda}(n,k)\binom{k+r}{k}k!x^{k}.\nonumber
\end{align}
Therefore, by \eqref{41}, we obtain the following theorem.
\begin{theorem}
For $n, r\ge 0$, we have
\begin{displaymath}
D_{r}F_{n,\lambda}(x)= \sum_{k=0}^{n}\binom{k+r}{k}k! S_{2,\lambda}(n,k) x^{k}.
\end{displaymath}
\end{theorem}
Let $p\in\mathbb{N}$. On the one hand, we have
\begin{align}
\quad \sum_{n=0}^{\infty}\binom{n+r}{r}(n)_{p,\lambda}x^{n}&=\sum_{n=0}^{\infty}\binom{n+r}{r}x^{n}\sum_{k=0}^{p}S_{2,\lambda}(p,k)(n)_{k}\label{42} \\
&=\sum_{k=0}^{p}S_{2,\lambda}(p,k)k!\sum_{n=k}^{\infty}\binom{n+r}{r}\binom{n}{k}x^{n}\nonumber\\
&=\sum_{k=0}^{p}S_{2,\lambda}(p,k)k!\sum_{n=k}^{\infty}\binom{k+r}{k}\binom{n+r}{n-k}x^{n} \nonumber \\
&=\sum_{k=0}^{p}S_{2,\lambda}(p,k)k!\binom{k+r}{r}x^{k}\sum_{n=0}^{\infty}\binom{n+k+r}{n}x^{n}\nonumber \\
&= \sum_{k=0}^{p}S_{2,\lambda}(p,k)k!\binom{k+r}{r}x^{k}\Big(\frac{1}{1-x}\Big)^{r+k+1} \nonumber \\
&=\Big(\frac{1}{1-x}\Big)^{r+1}D_{r}F_{p,\lambda}\Big(\frac{x}{1-x}\Big). \nonumber
\end{align}
On the other hand, we also have
\begin{align}
\Big(x\frac{d}{dx}\Big)_{p,\lambda}\Big(\frac{1}{1-x}\Big)^{r+1}&=\sum_{n=0}^{\infty}\binom{n+r}{n}\Big(x\frac{d}{dx}\Big)_{p,\lambda}x^{n}=\sum_{n=0}^{\infty}\binom{n+r}{n}(n)_{p,\lambda}x^{n}.\label{43}\\
\end{align}
Therefore, by \eqref{42} and \eqref{43}, we obtain the following theorem.
\begin{theorem}
For $p\in\mathbb{N}$, and $r\in\mathbb{Z}$ with $r\ge 0$, we have
\begin{displaymath}
\Big(x\frac{d}{dx}\Big)_{p,\lambda}\Big(\frac{1}{1-x}\Big)^{r+1}=\sum_{n=0}^{\infty}\binom{n+r}{n}(n)_{p,\lambda}x^{n}=\Big(\frac{1}{1-x}\Big)^{r+1}D_{r}F_{p,\lambda}\Big(\frac{x}{1-x}\Big).
\end{displaymath}	
\end{theorem}
By using Theorem 8, we observe that
\begin{align}
\frac{1}{r!}\Big(\frac{d}{dx}\Big)^{r}\bigg[\frac{x^{r}}{(1-x)^{2}}F_{p,\lambda}\Big(\frac{x}{1-x}\Big)\bigg]&= \frac{1}{r!}\Big(\frac{d}{dx}\Big)^{r}\bigg[\sum_{n=1}^{\infty}\Big((1)_{p,\lambda}+(2)_{p,\lambda}+\cdots+(n)_{p,\lambda}\Big)\bigg]x^{n+r}\label{44}\\
&=	\sum_{n=1}^{\infty}\Big((1)_{p,\lambda}+(2)_{p,\lambda}+\cdots+(n)_{p,\lambda}\Big)\binom{n+r}{r}x^{n}.\nonumber
\end{align}
Now, from Theorem 11 and \eqref{44} we obtain
\begin{align}
&\sum_{k=1}^{\infty}(k)_{p,\lambda}\bigg(\Big(\frac{1}{1-x}\Big)^{r+1}-\sum_{l=0}^{k}\binom{r+l}{l}x^{l}\bigg)\label{45}\\
&=\sum_{k=1}^{\infty}(k)_{p,\lambda}\sum_{n=k+1}^{\infty}\binom{n+r}{n}x^{n} =\sum_{n=1}^{\infty}\binom{n+r}{n}x^{n}\sum_{k=0}^{n-1}(k)_{p,\lambda} \nonumber \\
&=\sum_{n=1}^{\infty}\binom{n+r}{n}x^{n}\bigg(\sum_{k=1}^{n}(k)_{p,\lambda}-(n)_{p,\lambda}\bigg)\nonumber \\
&=\sum_{n=1}^{\infty}\binom{n+r}{r}\Big((1)_{p,\lambda}+(2)_{p,\lambda}+\cdots+(n)_{p,\lambda}\Big)x^{n}-\sum_{n=1}^{\infty}\binom{n+r}{r}(n)_{p,\lambda}x^{n}\nonumber \\
&=\frac{1}{r!}\Big(\frac{d}{dx}\Big)^{r}\bigg[\frac{x^{r}}{(1-x)^{2}}F_{p,\lambda}\Big(\frac{x}{1-x}\Big)\bigg]-\Big(x\frac{d}{dx}\Big)_{p,\lambda}\Big(\frac{1}{1-x}\Big)^{r+1}. \nonumber
\end{align}
Therefore, by \eqref{45}, we obtain the following differential equation for the degenerate Fubini polynomials.
\begin{theorem}
For $r\in\mathbb{Z}$ with $r\ge 0$, and $p\in\mathbb{N}$, we have 	
\begin{align*}
& \frac{1}{r!}\Big(\frac{d}{dx}\Big)^{r}\bigg[\frac{x^{r}}{(1-x)^{2}}F_{p,\lambda}\Big(\frac{x}{1-x}\Big)\bigg]-\Big(x\frac{d}{dx}\Big)_{p,\lambda}\Big(\frac{1}{1-x}\Big)^{r+1}\\
&=	\sum_{k=1}^{\infty}(k)_{p,\lambda}\bigg(\frac{1}{(1-x)^{r+1}}-\sum_{j=0}^{k}\binom{r+j}{j}x^{j}\bigg).
\end{align*}
\end{theorem}

\section{conclusion}
In this paper, we studied the degenerate Fubini and the degenerate Bell polynomials and found some further properties and identities on such polynomials which had been introduced by Kim-Kim. Especially, we found several expressions for generating functions for the sum $S_{n,\lambda}(p)=(1)_{p,\lambda}+(2)_{p,\lambda}+\cdots+(n)_{p,\lambda}$ of the values of the generalized falling factorials at positive consecutive integers. \par
As we mentioned earlier, there are various ways of studying special numbers and polynomials, including generating functions, $p$-adic analysis, umbral calculus, combinatorial methods, differential equations, special functions, probability theory, analytic number theory and operator theory. \par
It is one of our future projects to continue to explore various degenerate versions of many special polynomials and numbers by making use of such various tools.

\end{document}